\documentclass[12pt]{amsart}
\usepackage{amsmath}
\usepackage{amsfonts}
\usepackage{amssymb}
\usepackage{epsfig}
\newtheorem{lm}{Lemma}[section]
\newtheorem{tr}{Theorem}
\newtheorem{cor}{Corollary}[section]
\newtheorem{corr}{Corollary}

\newtheorem{pr}{Proposition}[section]
\usepackage{color}

\begin{document}

\title[Generators of supersymmetric polynomials]{Generators of supersymmetric polynomials in positive characteristic}
\author{A.N. Grishkov}
\address{Departamento de Matematica, Universidade de Sao Paulo, Caixa Postal 66281, 05315-970 - São Paulo, Brazil}
\email{shuragri@gmail.com}
\author{F. Marko}
\address{Penn State Hazleton, 76 University Drive, Hazleton PA 18202, USA}
\email{fxm13@psu.edu}
\author{A.N. Zubkov}
\address{Omsk State Pedagogical University, Chair of Geometry, 644099 Omsk-99, Tuhachevskogo Embankment 14, Russia}
\email{zubkov@iitam.omsk.net.ru}
\begin{abstract}
In \cite{kt}, Kantor and Trishin described the algebra of polynomial invariants of the adjoint representation of the Lie supergalgebra $gl(m|n)$ and a related algebra $A_s$ of what they called pseudosymmetric polynomials over an algebraically closed field $K$ of characteristic zero.
The algebra $A_s$ was investigated earlier by Stembridge who in \cite{s} called the elements of $A_s$ supersymmetric polynomials and determined  generators of $A_s$.

The case of positive characteristic $p$ has been recently investigated by La Scala and Zubkov in \cite{zs2}. They formulated two conjectures describing generators of polynomial invariants of the adjoint action of the general linear supergroup $GL(m|n)$ and generators of $A_s$, respectively.
In the present paper we prove both conjectures.
\end{abstract}
\maketitle

\section*{Introduction and notation}

Let $K$ be an algebraically closed field $K$ of positive characteristic $p\neq 2$.
The following notation is related to general linear supergroup $G=GL(m|n)$.
Let $K[c_{ij}]$ be a commutative superalgebra freely generated by elements $c_{ij}$ for $1\leq i,j \leq m+n$, where $c_{ij}$ is even if either $1\leq i,j\leq m$ or $m+1\leq i,j \leq m+n$, and $c_{ij}$ is odd otherwise.
Denote by $C$ the generic matrix $(c_{ij})_{1\leq i,j \leq m+n}$ and write $C$ as a block matrix $\begin{pmatrix}C_{00}&C_{01}\\C_{10}&C_{11}\end{pmatrix}$, where entries of $C_{00}$ and $C_{11}$ are even and entries of $C_{01}$ and $C_{10}$ are odd. The localization of $K[c_{ij}]$ by elements $det(C_{00})$ and $det(C_{11})$ is the coordinate superalgebra $K[G]$
of the general linear supergroup $G=GL(m|n)$.
The general linear supergroup $G=GL(m|n)$ is a group functor from the category $SAlg_K$ of
commutative superalgebras over $K$ to the category of groups, represented by
its coordinate ring $K[G]$, that is $G(A)=Hom_{SAlg_K}(K[G], A)$ for $A\in SAlg_K$.
Here, for $g\in G(A)$ and $f\in K[G]$ we define $f(g)=g(f)$.
Denote by $Ber(C)=det(C_{00} -C_{01}C_{11}^{-1}C_{10}) det(C_{11})^{-1}$ the Berezinian element.

Let $T$ be the standard maximal torus in $G$ and $X(T)$ be a set of characters.
Let $V$ be a $G$-supermodule with weight decomposition $V=\oplus_{\lambda\in X(T)} V_{\lambda}$, where $\lambda=(\lambda_1, \ldots, \lambda_{m+n})$ and
each $V_\lambda$ splits into a sum of its even subspace $(V_\lambda)_0$ and odd subspace $(V_\lambda)_1$.
The (formal) supercharacter $\chi_{sup}(V)$ of $V$ is defined as
\[\chi_{sup}(V)=\sum_{\lambda\in X(T)} (dim (V_\lambda)_0 - dim (V_\lambda)_1) x_1^{\lambda_1} \ldots x_m^{\lambda_m} y_1^{\lambda_{m+1}} \ldots y_n^{\lambda_{m+n}}.\]

If $V$ is a $G$-supermodule with a homogeneous basis $\{v_1,
\ldots, v_a, v_{a+1},$ $ \ldots , v_{a+b}\}$ such that $v_i$ is
even for $1\leq i\leq a$ and $v_i$ is odd for $a+1\leq i\leq a+b$,
and the image $\rho_V(v_i)$ of a basis element $v_i$ under a
comultiplication $\rho_V$ is given as $\rho_V(v_i)=\sum_{1\leq
j\leq a+b} v_j \otimes f_{ji}$, then the supertrace $Tr(V)$ is
defined as $\sum_{1\leq i\leq a} f_{ii} - \sum_{a+1\leq i\leq a+b}
f_{ii}$.

Let $E$ be a standard $G$-supermodule given by basis elements
$e_1, \ldots,$ $ e_m$ that are even and $e_{m+1}, \ldots, e_{m+n}$
that are odd and by comultiplication $\rho_E(e_i)=\sum_{1\leq
j\leq m+n} e_j \otimes c_{ji}$. Denote by $\Lambda^r(E)$ the
$r$-th exterior power of $E$ and by $C_r$ the supertrace of
$\Lambda^r(E)$.

The algebra $R$ of invariants with respect to adjoint action of $G$ is a set of functions $f\in K[G]$ satisfying $f(g_1^{-1}g_2g_1)=f(g_2)$ for any $g_1,g_2\in G(A)$ and any commutative supergalgebra $A$ over $K$. The algebra $R_{pol}$ of polynomial invariants is a subalgebra of $R$ consisting of polynomial functions.

As in the characteristic zero case, the description of polynomial invariants $R_{pol}$ can be reduced to that of the algebra $A_s$ of supersymmetric polynomials. The algebra $A_s$ consists of polynomials $f(x|y)=f(x_1, \ldots, x_m,y_1, \ldots y_n)$ that are symmetric in variables $x_1, \ldots x_m$
and $y_1, \ldots, y_n$ separately and such that $\frac{d}{dT}f(x|y)_{x_1=y_1=T}=0$.
The main tool for this reduction is the Chevalley epimorphism
$\phi: K[G]\to A$, where $A=K[x_1^{\pm 1}, \ldots, x_m^{\pm 1},y_1^{\pm 1} \ldots y_n^{\pm 1}]$, and $\phi(c_{ij})=\delta_{ij} x_i$ for $1\leq i\leq m$
and $\phi(c_{ij})=\delta_{ij} y_{i-m}$ for $m+1\leq i\leq m+n$.
Then for any $G$-supermodule $V$ we have $\phi(Tr(V))=\chi_{sup}(V)$.
In particular, for $0\leq r$ we have
\[\phi(C_r)=c_r=\sum_{0\leq i \leq min(r,m)} (-1)^{r-i} \sigma_i(x_1, \ldots, x_m)p_{r-i}(y_1, \ldots, y_n),\]
where $\sigma_i$ is the $i$-th elementary symmetric function and $p_j$ is the $j$-th complete symmetric function.

A homogeneous polynomial $f(x|y)=\sum a_{\lambda} x_1^{\lambda_1} \ldots x_m^{\lambda_m} y_1^{\lambda_{m+1}}\ldots y_n^{\lambda_{m+n}}$
is called $p$-balanced if $p|(\lambda_i+\lambda_j)$ whenever $1\leq i \leq m<j\leq m+n$ and $a_{\lambda}\neq 0$.
Denote by $A_s(p)$ the subalgebra of $A_s$ generated by $p$-balanced polynomials.

The following theorem is the main result of this paper.

\begin{tr}\label{prva}(Conjecture 5.2 of \cite{zs2})
The algebra $A_s$ is generated over $A_s(p)$ by elements $c_r$ for $r\geq 0$.
\end{tr}

For a matrix $M$, denote by $\sigma_i(M)$ the $i$-th coefficient
of the characteristic polynomial of $M$. Then all elements $C_r$,
$\sigma_i(C_{00})^p$, $\sigma_j(C_{11})^p$,
$\sigma_n(C_{11})^pBer(C)^k\in R_{pol}$. As a consequence of
previous theorem we obtain that these are generators of algebra
$R_{pol}$.

\begin{tr}(Conjecture 5.1 of \cite{zs2})
The algebra $R_{pol}$ is generated by elements
\[C_r,\sigma_i(C_{00})^p,\sigma_j(C_{11})^p,\sigma_n(C_{11})^pBer(C)^k,\]
where $0\leq r, 1\leq i\leq m, 1\leq j\leq n, 0<k<p$.
\end{tr}
\begin{proof}
Theorem 5.2 of \cite{zs2} states that the restriction of $\phi$ on
$R$ is a monomorphism. It was also noticed there that
$\phi(R_{pol})\subset A_s$ which is a consequence of arguments
analogous to Theorem 1.1 of \cite{kt}. Proposition 5.1 of
\cite{zs2} states that $\sigma_i(x_1, \ldots, x_m)^p$ for $1\leq
i\leq m$, $\sigma_j(y_1, \ldots, y_n)^p$ for $1\leq j\leq n$ and
$u_k(x|y)=\sigma_m(x_1, \ldots, x_m)^k\sigma_m(y_1, \ldots,
y_n)^{p-k}$ for $0<k<p$ are generators of algebra $A_s(p)$. Since
$\phi(C_r)=c_r$, $\phi(\sigma_i(C_{00})^p)=\sigma_i(x_1, \ldots,
x_m)^p$, $\phi(\sigma_j(C_{11})^p)=\sigma_j(y_1, \ldots, y_n)^p$
and $\phi(\sigma_n(C_{11})^pBer(C)^k)$ $=u_k(x|y)$, the statement
follows from Theorem \ref{prva}.
\end{proof}

\begin{corr}
Algebra $R$ is equaled to $R_{pol}[\sigma_m(C_{00})^{\pm p}, \sigma_n(C_{11})^{\pm p}]$.
\end{corr}
\begin{proof}
If $f\in R$, then its multiple by a sufficiently large power of
$\sigma_m(C_{00})^p \sigma_n(C_{11})^p$ is a polynomial invariant.
\end{proof}

\section{Nice supersymmetric polynomials}

In this section we will compare algebras corresponding to different values of $m,n$ and apply the Schur functor. Therefore we adjust the notation slightly to reflect the dependence on $m,n$. For example, we will write $R(m|n)$ instead of $R$ and $A_s(m|n)$ instead of $A_s$.

Denote by $A_{ns}(m|n)$ the algebra of "nice" supersymmetric
polynomials. It is a subalgebra of $A_s(m|n)$ generated by
polynomials \[c_r(m|n)=\sum_{0\leq i \leq r} (-1)^{r-i}
\sigma_i(x_1, \ldots, x_m)p_{r-i}(y_1, \ldots, y_n),\]
\[\sigma_i(x,m)^p=\sigma_i(x_1,\ldots, x_m)^p,
\sigma_j(y,n)^p=\sigma_j(y_1,\ldots, y_n)^p\] and
\[u_k(m|n)=\sigma_m(x,m)^k\sigma_n(y,n)^{p-k}\] for $1\leq i\leq m$,
$1\leq j\leq n$ and $0<k<p$. Since $\phi(C_r(m|n))=c_r(m|n)$,
using Proposition 5.1 of \cite{zs2} we conclude that
$A_{ns}(m|n)\subset \phi(R_{pol}(m|n))$. We will see below that
this inclusion is actually an equality. Denote by $A_{ns}(m|n,t)$
the homogeneous component of $A_{ns}(m|n)$ of degree $t$.

For any integers $M\geq m, N\geq n$ there is a graded superalgebra
morphism $p_e : K[x_1,\ldots,x_M, y_1, \ldots, y_N]\to K[x_1,\ldots, x_m, y_1,\ldots, y_n]$
that maps $x_i\mapsto x_i$ for $i\leq m$,
$y_j\mapsto y_j$ for $j\leq n$ and the remaining generators $x_i$, $y_j$ to zero.
Clearly $p_e$ restricts to a map from $A_s(M|N)$ into $A_s(m|n)$.

\begin{lm}\label{jedna}
The map $p_e$ takes $A_{ns}(M|N)$ to $A_{ns}(m|n)$.
\end{lm}
\begin{proof}
Assume $M>m$ or $N>n$. It follows from
$p_e(\sigma_i(x,M))=\sigma_i(x,m)$ if $i\leq m$ and
$p_e(\sigma_i(x,M))=0$ otherwise;
$p_e(\sigma_j(y,N))=\sigma_j(y,n)$ if $j\leq n$ and
$p_e(\sigma_j(y,N))=0$ otherwise; $p_e(c_r(M,N))=c_r(m|n)$ if
$r\leq m,n$ and $p_e(c_r(M,N))=0$ otherwise, and $p_e(u_k(M|N))$
$=0$.
\end{proof}

For the integers $M\geq m, N\geq n$ consider the Schur
superalgebra $S(M|N, r)$ and its idempotent
$e=\sum_{\mu}\xi_{\mu}$, where the sum runs over all weights $\mu$
for which $\mu_i=0$ whenever $m < i\leq M$ or $M+n < i\leq M+N$.
Then $S(m|n, r)\simeq eS(M|N, r)e$ and there is a natural Schur functor
$S(M|N, r)-mod\to S(m|n, r)-mod$ given by
$V\mapsto eV$. If $V$ is a $S(M|N,r)$-supermodule, then $eV$ is a
supersubspace of $V$ and therefore, $eV$ has the canonical $S(m|n,r)$-supermodule
structure.

\begin{lm}\label{dva}
The map $p_e$ induces an epimorphism of graded algebras
$\phi(R_{pol}(M|N))\to \phi(R_{pol}(m|n))$.
\end{lm}
\begin{proof} Applying the Chevalley map $\phi$ to the collection of
simple polynomial $G$-supermodules $L$ and using Theorem 5.3 of
\cite{zs2} we obtain that the algebra $\phi(R_{pol})$ is spanned
by the supercharacters $\chi_{sup}(L)$. If $\lambda$ is a highest
weight of $L$, then $\chi_{sup}(L)$ is a homogeneous polynomial of
degree $r=|\lambda|=\sum_{1\leq i\leq m+n}\lambda_i$. By the
standard property of a Schur functor there is a simple $S(M|N,
r)$-supermodule $L'$ such that $eL'\simeq L$. Since
$p_e(\chi_{sup}(L'))=\chi_{sup}(L)$, the claim follows.
\end{proof}

\begin{pr}$\phi(R_{pol}(m|n))=A_{ns}(m|n)$.
\end{pr}
\begin{proof} Fix a homogeneous element $f\in\phi(R_{pol}(m|n))$ of
degree $r$ and choose $M\geq m$ strictly greater than $r$.
By Lemma \ref{dva} there is a homogeneous polynomial $f'\in \phi(R_{pol}(M|n))$ of
degree $r$ such that $p_e(f')=f$.
Using Chevalley map and applying Theorem 5.3 of \cite{zs2} to the collection of costandard polynomial
modules $\nabla(\mu)$ we obtain that $f'$ is a linear combination of supercharacters $\chi_{sup}(\nabla(\mu))$, where
$\mu$ runs over polynomial dominant weights with $|\mu|=r$.

The assumption $M>r$, Theorem 5.4 and Proposition 5.6 of \cite{zs1} imply that for the highest weight $\mu=(\mu_+|\mu_-)$
we have $\mu_-=p\overline{\mu}$ and
$\nabla(\mu)\simeq\nabla(\mu_+|0)\otimes F(\overline{\nabla}(\overline{\mu}))$,
where $F$ is the Frobenius map and
$\overline{\nabla}(\overline{\mu})$ is the costandard $GL(n)$-module with the highest weight $\overline{\mu}$.
Therefore
\[\chi_{sup}(\nabla(\mu))=\chi_{sup}(\nabla(\mu_+|0))\chi(\overline{\nabla}(\overline{\mu}))^p.\]
Since $\chi(\overline{\nabla}(\overline{\mu}))^p$ is a polynomial in $\sigma_j(y,n)^p$, all that remains to show is that
$\chi_{sup}(\nabla(\mu_+|0))$ belongs to $A_{ns}(M|n)$ and use Lemma \ref{jedna}.

By Theorem 6.6 of \cite{zs1} the character $\chi_{sup}(\nabla(\pi|0))$ for a polynomial weight $\pi$ does not depend on
the characteristic of the ground field. Therefore we can temporarily assume that
$char K=0$. In this case the category $S(M|n, r)-mod$ is semisimple and its simple modules are
$\nabla(\lambda)$, where $\lambda$ runs over $(M|n)$-hook weights.

An exterior power $\Lambda^t(E(M|n))$ for $t\leq M$ has a unique maximal weight $(1^t |0)$.
Consequently, a $S(M|n, r)$-supermodule $V=$
\[\Lambda^M(E(M|n))^{\otimes\pi_M}\otimes\Lambda^{M-1}(E(M|n))^{\otimes
(\pi_{M-1}-\pi_M)}\otimes\ldots\otimes\Lambda^1(E(M|n))^{\otimes(\pi_1-\pi_2)}\]
has the unique maximal weight $(\pi|0)$ and supercharacter
\[\chi_{sup}(V)=c_1^{\pi_1-\pi_2}\ldots c_{M-1}^{\pi_{M-1}-\pi_M}c_M^{\pi_M}.\]

The module $V$ is a direct sum of $L(\pi|0)=\nabla(\pi|0)$ and $L(\kappa)=\nabla(\kappa)$ with $\kappa <(\pi|0)$.
Since $\kappa$ is a polynomial weight, it implies that $\kappa=(\kappa_+|0)$ and $\kappa_+<\pi$.
Using induction on $\pi$ we derive that $\chi_{sup}(\nabla(\pi|0))$ is a polynomial in $c_1,\ldots , c_M$
hence it belongs to $A_{ns}(M|n)$.
\end{proof}

\begin{cor}\label{Dos}
The morphism $p_e$ maps $A_{ns}(M|N,t)$ onto $A_{ns}(m|n,t)$.
\end{cor}

\section{Proof of Theorem \ref{prva}}

We will need the following crucial observation.

\begin{lm}\label{grish}
If $f\in A_s(m|n)$ is divided by $x_m$, then $f$ is divided by a nonconstant element of $A_{ns}(m|n)$.
\end{lm}
\begin{proof}
We can assume $f\neq 0$ and use the symmetricity of $f$ in variables $x_1,\ldots, x_m$ and $y_1, \ldots, y_n$ to
write $f=x_1^a \ldots x_m^a y_1^b \ldots y_n^b g$, where exponents $a>0,b\geq 0$ and a polynomial $g$ such that $g|_{x_m=y_n=0}\neq 0$ are unique.
Then
$$f|_{x_m=y_n=T}=T^{a+b}x_1^a \ldots x_{m-1}^a y_1^b \ldots y_{n-1}^b g|_{x_m=y_n=T}$$
$$= T^{a+b}x_1^a \ldots x_{m-1}^a y_1^b \ldots y_{n-1}^b g_0
+ T^{a+b+1}x_1^a \ldots x_{m-1}^a y_1^b \ldots y_{n-1}^b g_1,$$
where we write $g|_{x_m=y_n=T}=g_0+Tg_1$.
The requirement $g|_{x_m=y_n=0}\neq 0$ implies $g_0\neq 0$. Since $\frac{d}{dT}f|_{x_m=y_n=T}=0$, this is only possible if
$a+b \equiv 0 \pmod p$. Since $a>0$, the polynomial $x_1^a\ldots x_m^ay_1^b\ldots y_n^b$ is not constant, and is a product of
$\sigma_m(x,m)^p$, $\sigma_n(y,n)^p$ and $u_k(m|n)$ which belongs to $A_{ns}(m|n)$.
In fact, since $a>0$, we have that $f$ is divisible either by $\sigma_m(x,m)^p$ or by some $u_k(m|n)$.
\end{proof}

{\it Proof of Theorem \ref{prva}.}
Using Proposition 5.1 of \cite{zs2} the statement of the theorem is equivalent to the equality $A_s(m|n)=A_{ns}(m|n)$.

Fix $n$ and assume that $m$ is minimal such that there exists a polynomial $f\in A_s(m|n)\setminus A_{ns}(m|n)$ and choose $f$ such that it is homogeneous and of the minimal degree.
Then its reduction $f|_{x_m=0}\in A_{ns}(m-1|n)$ is a nonzero polynomial $h(c_t(m-1|n),\sigma_i(x,m-1)^p, \sigma_j(y,n)^p, u_k(m-1|n))$ in elements $c_t(m-1|n)$, $\sigma_i(x,m-1)^p$, $\sigma_j(y,n)^p$ and $u_k(m-1,n)$ where $t\geq 0$, $1\leq i\leq m-1$, $1\leq j\leq n$ and $0<k<p$.
By Corollary \ref{Dos} there are elements $v_k\in A_{ns}(m|n)$ of degree $mk+(p-k)n$ such that $v_k|_{x_m=0}=u_k(m-1|n)$. Since $c_t(m|n)|_{x_m=0}=c_t(m-1|n)$,
$\sigma_i(x,m)^p|_{x_m=0}=\sigma_i(x,m-1)^p$ and $\sigma_j(y,n)^p|_{x_m=0}=\sigma_j(y,n)^p$, the polynomial
$l=f-h(c_t(m|n), \sigma_i(x,m)^p, \sigma_j(y,n)^p, v_k(m|n))$ satisfies $l|_{x_m=0}=0$.
Since the degree of $l$ does not exceed the degree of $f$, $l\in A_s(m|n)$ and $x_m$ divides $l$, Lemma \ref{grish} implies
that $l=l_0l_1$, where $l_0\in A_{ns}(m|n)$ and the degree of $l_1$ is strictly less than the degree of $f$.
But $l_1\in A_s(m|n)\setminus A_{ns}(m|n)$ which is a contradiction with our choice of $f$.
\qed

\section{Elementary proof of Theorem \ref{prva}}

A closer look at the proof of Theorem \ref{prva} reveals that Corollary \ref{Dos} is the only result from Section 1 that was used in it.
Actually, only the following weaker statement was required in the proof of Theorem \ref{prva}.

\begin{pr}\label{elem}
For each $0<k<p$ there is a polynomial $v_k\in A_{ns}(m|n)$ of degree $(m-1)k+(p-k)n$ such that $v_k|_{x_m=0}=u_k(m-1|n)$.
\end{pr}

In this section we give a constructive elementary proof of Proposition \ref{elem} that bypasses a use of the Schur functor and results about costandard modules derived in \cite{zs1}.

Fix $0<k<p$ and denote $s=\lceil \frac{k}{p-k}\rceil$.  Then for $i=0, \ldots s-1$ define $k_i=(i+1)k-ip>0$
and $k_p=sp-(s+1)k\geq 0.$ The relations
\[k_i+(p-k)=k_{i-1},\, k_p+k=s(p-k),\,k_i+k_p=(s-i)(p-k)\]
will be used repeatedly.

A symbol $\Delta$ will denote a nondecreasing sequence $(i_1\leq \ldots \leq i_t)$ of natural numbers, where $0\leq t<s$. We denote $||\Delta||=t$ and $|\Delta|=\sum_{j=1}^ti_j$. In particular, we allow $\Delta=(\emptyset)$ and set $||(\emptyset)||=|(\emptyset)|=0$.
Denote by $Supp(\Delta)$ a maximal increasing subsequence of $\Delta$.
Further, denote by $Sym_M(x_1^{a_1}\dots x_r^{a_r})$ and $Sym_N(y_1^{b_1}\dots y_r^{b_r})$ respectively
a homogeneous symmetric polynomial in variables  $x_1,\ldots, x_M$  and $y_1, \ldots, y_N$ respectively with a general monomial term
$x_1^{a_1}\dots x_r^{a_r}$ and $y_1^{b_1}\dots y_r^{b_r}$ respectively.
Denote
\[\begin{array}{l}
(\Delta,j)_{M,N}=\\
Sym_M(x_1^k\dots x_{M-t}^kx_{M-t+1}^{k_{i_1}}\dots x_{M}^{k_{i_t}})Sym_N(y_1^{p-k}\dots y_{N-j-1}^{p-k}y_{N-j}^{k_{p}})
\end{array}\]
for $0\leq ||\Delta||=t \leq M$ and $0\leq j<N$, and $(\Delta,l)_{M,N}=0$ otherwise;
\[\begin{array}{l}
[\Delta,j]_{M,N}=\\
Sym_M(x_1^k\dots x_{M-t}^kx_{M-t+1}^{k_{i_1}}\dots x_{M}^{k_{i_t}})Sym_N(y_1^{p-k}\dots y_{N-j}^{p-k})
\end{array}\]
for $0\leq ||\Delta||=t \leq M$ and $0\leq j \leq N$, and $[\Delta,j]_{M,N}=0$ otherwise;
\[\begin{array}{l}
\{\Delta,l,j\}_{M,N}=\\
Sym_M(x_1^k\dots x_{M-t-1}^kx_{M-t}^{l(p-k)}x_{M-t+1}^{k_{i_1}}\dots x_{M}^{k_{i_t}})Sym_N(y_1^{p-k}\dots y_{N-j}^{p-k})
\end{array}\]
for $0\leq ||\Delta||=t < M$ and $0\leq j \leq N$  and any $l$, and $\{\Delta,l,j\}_{M,N}=0$ otherwise.

For $f \in K[x_1, \ldots, x_m,y_1, \ldots, y_n]$ define $\psi(f)=f|_{x_m=y_n=T}$ and
for $g, h\in K[x_1, \ldots, x_{m-1},y_1, \ldots, y_{n-1},T]$ write $g \equiv h$ if and only if $\frac{d}{dT}(g-h)=0.$

For simplicity write $(\Delta,j)$, $[\Delta,j]$ and
$\{\Delta,l,j\}$ short for $(\Delta,j)_{m-1,n-1}$,
$[\Delta,j]_{m-1,n-1}$ and $\{\Delta,l,j\}_{m-1,n-1}$.

\begin{lm}\label{l1} We have

$\begin{array}{l}
\psi\{\Delta,l,j\}_{m,n}\equiv \\
T^k\{\Delta,l,j-1\}+T^{(l+1)(p-k)}[\Delta,j]+T^{l(p-k)}[\Delta,j-1]+\\
\sum_{i\in Supp(\Delta)}(T^{k_i}\{\Delta\setminus i,l,j-1\}+T^{k_{i-1}}\{\Delta\setminus i,l,j\}).
\end{array}$

and

$\begin{array}{l}
\psi(\Delta,j)_{m,n}\equiv \\
T^k(\Delta,j-1)+T^{s(p-k)}[\Delta,j]+\\
\sum_{i\in Supp(\Delta)}(T^{k_i}(\Delta\setminus i,j-1)+T^{k_{i-1}}(\Delta\setminus i,j)+T^{(s-i)(p-k)}[\Delta\setminus i,j]).
\end{array}$
\end{lm}
\begin{proof}
The first relation follows from

$\begin{array}{l}
\psi(Sym_m(x_1^k\dots x_{m-t-1}^kx_{m-t}^{l(p-k)}x_{m-t+1}^{k_{i_1}}\dots x_{m}^{k_{i_t}}))=\\
\delta_{t,m-1}T^k Sym_{m-1}(x_1^k\dots x_{m-1-t-1}^kx_{m-1-t}^{l(p-k)}x_{m-1-t+1}^{k_{i_1}}\dots x_{m-1}^{k_{i_t}})\\
+T^{l(p-k)}Sym_{m-1}(x_1^k\dots x_{m-t-1}^k x_{m-t}^{k_{i_1}}\dots x_{m-1}^{k_{i_t}})+\\
\sum_{i\in Supp(\Delta)}T^{k_i}Sym_{m-1}(x_1^k\dots
x_{m-t-1}^kx_{m-t}^{l(p-k)}x_{m-t+1}^{k_{i_1}}\dots
\widehat{x_{m-t+i}^{k_i}} \dots x_{m-1}^{k_{i_t}}) ,\end{array}$

$\begin{array}{l}\label{tech}
\psi(Sym_n(y_1^{p-k}\dots y_{n-j}^{p-k}))=\\
\delta_{j,n}T^{p-k}Sym_{n-1}(y_1^{p-k}\dots y_{n-1-j}^{p-k})+\delta_{j,0}Sym_{n-1}(y_1^{p-k}\dots y_{n-j}^{p-k})
\end{array}$

and definitions of $(\Delta,j)$, $[\Delta,j]$ and $\{\Delta,l,j\}$.

Second relations follows from

$\begin{array}{l}
\psi(Sym_m(x_1^k\dots x_{m-t}^kx_{m-t+1}^{k_{i_1}}\dots x_{m}^{k_{i_t}}))=\\
\delta_{t,m}T^k Sym_{m-1}(x_1^k\dots x_{m-t-1}^kx_{m-t}^{k_{i_1}}\dots x_{m-1}^{k_{i_t}})\\
+\sum_{i\in Supp(\Delta)}T^{k_i}Sym_{m-1}(x_1^k\dots x_{m-t}^kx_{m-t+1}^{k_{i_1}}\dots \widehat{x_{m-t+i}^{k_i}} \dots x_{m-1}^{k_{i_t}})
,\end{array}$

$\begin{array}{l}
\psi(Sym_n(y_1^{p-k}\dots y_{n-j}^{p-k}y_{n-1-j}^{k_p}))=\\
\delta_{j,n-1}T^{p-k}Sym_{n-1}(y_1^{p-k}\dots y_{n-1-j-1}^{p-k}y_{n-1-j}^{k_p})\\
+T^{k_p}Sym_{n-1}(y_1^{p-k}\dots y_{n-1-j}^{p-k})
+\delta_{j,0}Sym_{n-1}(y_1^{p-k}\dots y_{n-1-j}^{p-k}y_{n-j}^{k_p})
\end{array}$

and definitions of $(\Delta,j)$, $[\Delta,j]$ and $\{\Delta,l,j\}$.
\end{proof}

Let us define

$\begin{array}{ll}
w=&\sum_{l=1}^{s-1}\sum_{l-n\leq |\Delta|\leq l}(-1)^{|\Delta|+s+l}(s-l)\{\Delta,l,l-|\Delta|\}_{m,n}\\
&+\sum_{s-1-n<|\Delta|<s}(-1)^{|\Delta|}(\Delta,s-1-|\Delta|)_{m,n}.
\end{array}$.

This expression is identical to

$\begin{array}{ll}
&\sum_{l=1}^{s-1}\sum_{0\leq |\Delta|\leq l}(-1)^{|\Delta|+s+l}(s-l)\{\Delta,l,l-|\Delta|\}_{m,n}\\
&+\sum_{\Delta}(-1)^{|\Delta|}(\Delta,s-1-|\Delta|)_{m,n}
\end{array}$

obtained by inserting additional terms that are all equal to zero.
Then

$\begin{array}{ll}
\psi(w)=&\sum_{l=1}^{s-1}\sum_{0\leq |\Delta|\leq l}(-1)^{|\Delta|+s+l}(s-l)\psi\{\Delta,l,l-|\Delta|\}_{m,n}\\
&+\sum_{\Delta}(-1)^{|\Delta|}\psi(\Delta,s-1-|\Delta|)_{m,n}
\end{array}$

which by Lemma \ref{l1} equals

$\begin{array}{l} \sum_{l=1}^{s-1}\sum_{0\leq |\Delta|\leq
l}(-1)^{s+l+|\Delta|}(s-l)
\Big(T^k\{\Delta,l,l-|\Delta|-1\}\\+T^{(l+1)(p-k)}[\Delta,l-|\Delta|]
+T^{l(p-k)}[\Delta,l-|\Delta|-1]\\
+\sum_{i\in Supp(\Delta)}(T^{k_i}\{\Delta\setminus i, l,l-|\Delta|-1\}+T^{k_{i-1}}\{\Delta\setminus i,l,l-|\Delta|\})\Big)\\

+\sum_{\Delta}(-1)^{|\Delta|} \\
\Big(T^k(\Delta,s-|\Delta|-2)+T^{s(p-k)}[\Delta,s-1-|\Delta|]\\
+\sum_{i\in Supp(\Delta)}(T^{k_i}(\Delta\setminus i,s-|\Delta|-2)+T^{k_{i-1}}(\Delta\setminus i,s-|\Delta|-1)\\
+T^{(s-i)(p-k)}[\Delta\setminus i, s-|\Delta|-1])\Big).
\end{array}$

To analyze this expression, denote by  $C(a)$ the coefficient corresponding to a term $a$ of $\psi(w)$.

Then $C(T^{k_i}(\Delta,s-2-i-|\Delta|)=(-1)^{|\Delta|+i} +(-1)^{|\Delta|+i+1}=0$ for each $i=0, \ldots, s-2$ and admissible $\Delta$ and
$C(T^{k_i}\{\Delta,l,l-1-|\Delta|-i\})=(-1)^{s+l+|\Delta|+i}+ (-1)^{s+l+|\Delta|+i+1}=0$ for each $i=0, \ldots, s-1$ and admissible $\Delta$.

Next, $C(T^{s(p-k)}[\Delta,s-|\Delta|-1])=(-1)^{|\Delta|+s+s-1}(s-(s-1))+(-1)^{|\Delta|}=0$ and if
$l-1-|\Delta|\neq 0$ or $\Delta\neq \emptyset$, then
$C(T^{l(p-k)}[\Delta,l-|\Delta|-1])=(-1)^{|\Delta|+s+l}(s-l)+(-1)^{|\Delta|+s+l-1}(s-l+1)+
(-1)^{|\Delta_l|},$ where $\Delta_l=\Delta\cup (s-l).$ Since $(-1)^{|\Delta_l|}=(-1)^{|\Delta|+s-l}$ we conclude that
$C(T^{l(p-k)}[\Delta,l-|\Delta|-1])=0$.
Finally, $C(T^{p-k}[\emptyset,0])=(-1)^{s+1}(s-1)+(-1)^{|\Delta_1|}=(-1)^{s+1}s.$
Therefore

$\begin{array}{ll}\psi(w)&=(-1)^{s+1}sT^{p-k}[\emptyset,0]\\
&=(-1)^{s+1}sT^{p-k}Sym_{m-1}(x_1^k\dots x_{m-1}^k)Sym_{n-1}(y_1^{p-k}\dots y_{n-1}^{p-k}).
\end{array}$

We can now easily prove Proposition \ref{elem}.

\medskip

{\it Proof of Proposition \ref{elem}.}
Since $s<p$ we can take
\[v_k=\frac{(-1)^{s}}{s}w+Sym_m(x_1^k\dots x_{m-1}^k)(y_1^{p-k}\dots y_n^{p-k}).\]
Then $\psi(v_k)=0$, hence $v_k\in A_s(m|n)$ and $v_k|_{x_m=0}=u_k(m-1|n)$.
It is easy to check that $v_k$ is homogeneous of degree $(m-1)k+(p-k)n$.
\qed.

\section{Concluding remarks}

The proof of Theorem \ref{prva} was influenced by Theorem 1 in \cite{s}. Although our proof uses different arguments we would like to remark that an analogue of Theorem 1 of \cite{s} remains valid over arbitrary commutative ring $A$ of any characteristic.

\begin{pr}\label{stem}
The algebra of polynomials $f(x_1, \ldots, x_m,y_1, \ldots, y_n)$ over a commutative ring $A$, symmetric in variables $x_1, \dots, x_m$ and
$y_1, \ldots, y_n$ separately and such that $f|_{x_m=y_n=T}$ does not depend on $T$ is generated by polynomials $c_r(x|y)$.
\end{pr}
\begin{proof}
Proof is a complete analogue of Theorem 1 of \cite{s}. Just replace every appearance of $\sigma_{m,n}^{(r)}$ by $c_r(x|y)$.
\end{proof}

Let us comment that if characteristic of $A$ is positive, then condition that $f|_{x_m=y_n=T}$ does not depend on $T$ is stronger than $\frac{d}{dT} f|_{x_m=y_n=T}=0$.

Proposition 3.1 of \cite{kt} states that, in the case of characteristic zero, the algebra $A_s$ is infinitely generated.
In the case of positive characteristic we have the following.

\begin{pr}
Algebra $A_s$ is finitely generated.
\end{pr}
\begin{proof}
The algebra $A_s$ is contained in $B=K[\sigma_i(x|m), \sigma_j(y|n)|1\leq i\leq m, 1\leq j\leq n]$.
Algebra $B$ is finitely generated over its subalgebra $B'=  K[\sigma_i(x|m)^p, \sigma_j(y|n)^p |1\leq i\leq m, 1\leq j\leq n]$ hence is a Noetherian $B'$-module.
However, $A_s$ contains $B'$ and is therefore finitely generated $B'$-module. Since $B'$ is finitely generated, so is $A_s$.
\end{proof}

\end{document}